\documentclass[11pt,english]{article}
\usepackage[T1]{fontenc}
\usepackage[latin9]{inputenc}
\usepackage{amsmath}
\usepackage{graphicx,pstricks}
\usepackage{amssymb}
\usepackage[arrow, matrix, curve]{xy}
\usepackage{mathabx} \usepackage{color} \usepackage{ifsym}
\usepackage[footnotesize]{caption} %cosi` sono 2pt in meno del testo, se vuoi solo 1pt devi mettere small; se vuoi allineare la didascalia in verticale all'altezza dei ":", invece che parta ad inizio riga quando vai a capo, nelle opzioni devi aggiungere hang

\makeatletter
\usepackage{epsfig}\usepackage{graphics}
\newtheorem{theorem}{Theorem}
\newtheorem{lem}{Lemma}[section]
\newtheorem{corollar}[lem]{Corollary}

\newcommand{\N}{{\mathbb N}}

\newcommand{\R}{{\mathbb R}}
\newcommand{\Z}{{\mathbb Z}}

\newcommand{\E}{{\mathbb{E}}}

\def\N{{\mathbb N}}
\def\Z{{\mathbb Z}}
\def\R{{\mathbb R}}
\def\P{{\mathbb P}}

\def\E{{\mathbb E}}

\def\0{{\bf 0}}

\def\dist{{\rm dist}}

\def\a{\alpha}

\def\d{\delta}
\def\e{\varepsilon}

\def\phi{\varphi}

\def\l{\lambda}

\def\s{\sigma}
\def\t{\tau}

\def\O{\Omega}

\def\T{\T}

\def\PP{{\cal P}}

\def\Cox{\hfill \Box}

\catcode`@=11 \@addtoreset{equation}{section} \catcode`@=12

\usepackage{babel}

\begin{document}

\title{Dependent particle deposition on a graph:\\  concentration properties of the height profile}

\author{S.R. Fleurke%
\thanks{Radiocommunications Agency Netherlands, Postbus 450, 9700 AL Groningen, The Netherlands,
\texttt{sjoert.fleurke@at-ez.nl} %
},\, 
M. Formentin\thanks{Ruhr-Universit\"at Bochum, Fakult\"at f\"ur Mathematik, Universit\"atsstrasse
150, 44780 Bochum, Germany, \texttt{Marco.Formentin@rub.de} %
} , \,and C. K\"ulske %\thanks{Research supported by Deutsche Forschungsgemeinschaft}
\thanks{Ruhr-Universit\"at Bochum, Fakult\"at f\"ur Mathematik, Universit\"atsstrasse
150, 44780 Bochum, Germany, \texttt{Christof.Kuelske@rub.de} %
} }
\maketitle
\begin{abstract} 
We present classes of models in which particles are dropped on an arbitrary fixed finite 
connected graph, obeying adhesion rules with screening. 
We prove that there is an invariant distribution for the resulting height profile, and Gaussian concentration 
for functions depending on the paths of the profiles. 
As a corollary we obtain a law of large numbers for the maximum height. {This describes the asymptotic speed with which the maximal height increases.} 

{The results incorporate the case of  independent particle droppings but extend to droppings according
to a driving Markov chain, and to droppings with possible deposition below the top 
layer up to a fixed finite depth, obeying a non-nullness condition for the screening rule.}  
The proof is based on an analysis of the Markov chain on height-profiles using coupling methods.
We construct a finite communicating set of configurations of profiles to which the chain 
keeps returning.  
\end{abstract}
\smallskip{}

%\noindent \vspace{10mm}
 \textit{Key--Words:} Random sequential adsorption, Particle deposition, Driven interfaces, 
 Particle systems, Gaussian bounds,  Concentration estimates, Coupling, Law of large numbers.

\section{Introduction}

Stochastic models for particle deposition have enjoyed much interest over the years,
motivated by applications ranging from car parking, physical chemistry to frequency assignment  \cite{Cha07, CoRe63, DeFl07, Evans,  Me83, Renyi, ShVo09}. 
In a series of papers particle deposition models with a number of different 
deposition rules were considered and exact solutions for models in solvable geometries were given 
\cite{Pr04, DeFlKu08,FlKu09,FlKu10,Gouet, PeSu05,Sudbury}. 
Natural probabilistic questions to be studied in cases where no closed solutions 
are available are limit laws for such processes in space \cite{PeYu02, ScPeYu07}   
or time. In particular one would like to have a law of large numbers for the maximum of the height variables 
and the behavior of the active or top region. 
Moreover there is a branch in probability which is interested  
in the investigation in concentration of measure properties for Markov chains 
and multidimensional stochastic processes \cite{CCKR06, ChRe09,Ku03,Ma98}, 
and we also want to look at deposition models in this spirit. 

In the present paper we consider models of discrete-time Markov chains describing 
the growth of adsorbed particles on a substrate.  
In our main example particles are dropped on the vertices 
of a finite connected graph $V$ according to a discrete time Markov chain
and obeying screening rules of adsorption.   
The particles pile up to integer heights according to an exclusion interaction  
between sites which are connected in $V$. 
Our last example softens the screening rule to allow adsorption below the top layer. 

We prove a strong law of large numbers for the maximal height and show convergence 
of the height profile to a stationary state. As the number of deposited particles grows linearly in time when we 
keep the graph (and hence the volume) fixed, we will look at the heights differences relative to the maximum. 
This map from height configurations to relative heights is just the same 
as the map from interface configurations to gradient configurations 
considered in models of interfaces in a Gibbs state \cite{EnKu08,FunSp, Vel06} 
when issues of stability of interfaces in the large volume limit are considered. 

Now, in our situation we show the convergence of the height-profile as seen from the maximum to an
invariant distribution using a coupling method.   
Our Markov chain has an unbounded state space, but the coupling turns 
out to be very good, namely we are able to show that the distribution of the coupling 
time can be controlled uniformly in the initial configurations. The physical reason for this 
is the following: however rough a profile is, there is always a chain of particle droppings 
which will make it flat and thereby erase the memory on the past. 
An essential ingredient for this to turn into a proof in the context of the general models 
we consider, is the construction of a finite set of profiles 
the chain communicates to in a time $s$ which is uniform in any starting configuration. 
The construction of this set is slightly subtle in the case of a non-i.i.d. chain of particle 
droppings where it is based on irreducibility and lazyness of the driving chain. 
In particular, from this  coupling the law of large numbers for the maximal height 
follows as a corollary from concentration results for path observables. 

\section{The models and the main results}

\subsection{Independent particle droppings}

Let $G=(V,E)$ be a finite connected graph. Write $i\sim j$ if $\{i,j\}\in E$, that is $i,j$ are adjacent. 
Consider the Markov chain
on the state space $\O:=\N_{0}^{V}$ of height configurations $h=(h_j)_{j\in V}$ 
obtained by choosing a site $i\in V$ according to
a probability $p(i)>0$, where $p\in\PP(V)$ is fixed, and adding a
particle at $i$ at height $\max \{h_j, \dist (j,i)\leq 1\}+1$ where $h_j$ is the maximum height 
at which a particle is already present at site $j$. 

The formal definition is as follows. Denote by 
$T_i : \O\rightarrow \O$ the operator which assigns to a configuration $h$ the configuration 
 $T_i h$ which is the configuration obtained by adding a particle at $i$, i.e. 
\begin{equation}
\begin{split}
(T_i h)_j&= \begin{cases}
 \max\{h_k: \dist (k,i)\leq 1   \}+1  & \text{if }  j= i \cr
 h_j & else
\end{cases}
\end{split}
\end{equation}

Look at the discrete time Markov chain with transition matrix $(M(h,h'))_{h,h'\in \O}$ given by 
\begin{equation}
\begin{split}
M(h,h')=\begin{cases}
p(i) & \text{if }  h' = T_i h \cr
0 & else
\end{cases}
\end{split}
\end{equation}
We denote the value of the configuration 
at time $t$ by $h(t)=(h_{i}(t))_{i\in V}$.

The model has the following property:
If $h'\in \O$ is such that  $h'_j = h_j + c$ for all $j\in V$ 
we have that $(T_i h')_j = (T_i h)_j + c$ and hence we can define the action of $T_i$ 
also on equivalence classes of height-profiles 
w.r.t. constant shifts $c$. { Let us extend the local state space to $\Z$ and allow for arbitrary $c\in\Z$.} We may choose then a representative 
of these equivalence classes in such a way that the height profile is zero at the maximum and 
negative elsewhere. That is, we introduce the variable $x_{i}=h_{i}-\max_{j\in V}h_{j}$. This is the
height profile seen from the maximum.

According to the exclusion rules
the process on $x=(x_{i})_{i\in V}$ is a Markov chain again, now with
state space $S:=(-\N_{0})^{V}$ and transition matrix 
$M(x,x')=M(h,h')$ when $x$ is the equivalence class of $h$ and 
$x'$ is the equivalence class of $h'$. 
We will show convergence to an invariant distribution of this Markov chain. In order to do this 
we need to prove recurrence, and therefore we need to make use of
the exclusion rules. To compare, consider the process in which particles
are added without exclusion. Then the distribution of the heights
becomes multinomial and the corresponding $x$-distribution won't
stabilize but have fluctuations of the order of the square-root of
the discrete time $n$.

\subsection{Markov chain particle droppings}
{Now the probability where to drop the next particle depends on where the last time a particle has fallen.}

Let $v(t)$ denote a Markov chain with state space $V$ and transition 
matrix $$A_{v,v'}=\P(v(t+1)=v'| v(t)=v).$$ We call this the driving Markov chain. 

We assume that $A=(A(v,v'))_{v,v'\in V}$ is {\em{irreducible}} (meaning that for all $v\neq v'$ there 
exists a time $s(v,v')$ such that $A^{s(v,v')}(v,v')>0$) and that it is {\em{lazy}} 
(meaning that $A(v,v)>0$ for all $v\in V$.)

This time look at the Markov chain $(h(t),v(t))$ with transition matrix 
$$(M(h,v;h',v'))_{h,v;h',v'\in \O\times V}$$ given by 
\begin{equation}
\begin{split}
M(h,v;h',v')=\begin{cases}
 A(v,v')& \text{ if }  h' = T_{v'} h  \cr
0 & \text{ else}
\end{cases}
\end{split}
\end{equation}
We denote the value of the configuration 
at time $t$ by $(h(t),v(t))$.

\subsection{Main results}

Our main goal will be the following theorem which provides a concentration 
estimate for a specific important example of an observable. Generalizations to other observables 
will become clear from the proof.

\begin{theorem} Assume that we are given either a model of independent particle droppings or, 
more generally Markov chain particle droppings on a connected graph with more than two vertices. 
Then the following holds. 

\begin{enumerate}
\item $x(t)$ converges in law to an invariant distribution, independently of the starting configuration. 

\item Define $m_V(t)=\max_{j\in V}h_{j}(t)$ to be the total height of the particle profile.  
Then there exists a positive constant $c$, depending on the model, such that 
\begin{equation}\label{bound1}
\begin{split} 
\P(|m_V(t) -\E m_V(t)| > y) &
\leq 2 \exp\Bigl( 
- \frac{c y^2 }{2 t} 
\Bigr) 
\end{split}
\end{equation}
where the bounds hold either if we take for $\P=\P_{\pi}$ the chain in equilibrium, or the chain 
started in any initial configuration.  

\item There exists a constant $C$ such that 
\begin{equation}\label{bound2}
\begin{split} 
&\sup_{t}\Bigl |\E_{\pi} m_V(t)- \E_0 m_V(t)\Bigr|\leq C\cr 
 \end{split}
\end{equation}
is uniformly bounded, where $\E_0 $ denotes the chain started in the flat configuration $h_j=0$ for all $j\in V$. 
\end{enumerate}
\end{theorem}
From the Theorem follows the SLLN for the variable $\frac{m_V(t)}{t}$ as $t$ tends 
to infinity and also the independence 
of the initial configuration.

 % Put also 
%$$I(h):= \text{argmax}_{i}h_i\subset V$$ for the set of maximizers 
%of the height function and similarly $$I(x):= \text{argmax}_{i}x_i\subset V$$ 

\section{Independent particle droppings - the proof }

We will now give a self-contained presentation of the proof for the  
first example of independent particle droppings. 

\subsection{Construction of communicating set - convergence to invariant distribution}

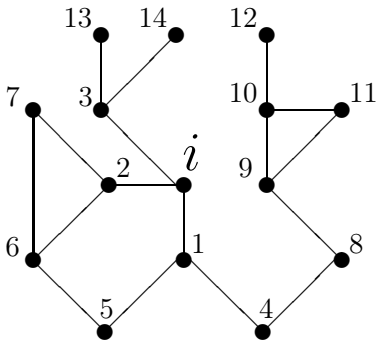
\begin{figure}
\setlength{\unitlength}{1mm}
\begin{center}
\begin{picture}(50,60)(0,-7)
% links
\multiput(5,10)(0,20){2}{\circle*{2}} %
\put(5,10){\line(1,1){10}} \put(5,10){\line(1,-1){10}}
\put(5,30){\line(0,-1){20}}
\put(5,30){\line(1,-1){10}} \put(15,20){\circle*{2}}
% midden
\put(24,20){\line(-1,0){10}} \put(24,20){\line(-1,1){10}}
\put(25,19){\line(0,-1){8}} \put(25,20){\circle*{2}}
% onder
\put(25,10){\circle*{2}} 
\put(24,10){\line(-1,-1){10}} %\put(25,19){\line(0,-1){9}}
\put(26,10){\line(1,-1){10}} \put(14.5,0.5){\circle*{2}}
\put(36,1){\line(1,1){10}}
\put(35.5,0.5){\circle*{2}}
% rechts
\put(36,20){\circle*{2}} \put(36,20){\line(1,1){10}}
\put(36,20){\line(1,-1){10}} \put(46,30){\circle*{2}}
\put(46,10){\circle*{2}} \put(36,30){\line(1,0){10}}
\put(36,30){\line(0,-1){10}}
% boven
\put(36,30){\circle*{2}} \put(36,30){\line(0,1){10}}
\put(36,40){\circle*{2}} %
\put(14,30){\circle*{2}} \put(14,30){\line(0,1){10}}
\put(14,30){\line(1,1){10}} \put(14,40){\circle*{2}}
\put(24,40){\circle*{2}} %\put(24,40){\line(1,0){12}}
% tekst
\put(24.8,22){\huge{$i$}} \put(26,11){1} \put(13.8,2.3){5}
\put(35,2.3){4}
\put(1.3,10.5){6}
\put(1.3,30.5){7}
\put(16,21){2}
\put(11,30.5){3}
\put(9,41){13}
\put(19,41){14}
\put(47,11){8}
\put(32.3,21){9}
\put(31,31){10}
\put(31,41){12}
\put(47,31){11}
\end{picture}
\caption{An example of an $i$-ordering.} \label{fig:1}
\end{center}
\end{figure}

For each vertex $i\in V$ we pick an $i$-dependent 
ordering $a^{(i)}=(a^{(i)}_{1},\dots,a^{(i)}_{|V|-1})$
of the sites in the set $V\backslash \{i\}$, 
starting with $a^{(i)}_{1}$ to be a nearest neighbor of $i$ 
and the additional property that $d(a^{(i)}_{k},\{a^{(i)}_{1},\dots,a^{(i)}_{k-1}\})=1$
where $d$ is the graph distance (see fig. \ref{fig:1}). 
This means that $a^{(i)}$ describes a way how 
the set $V$ can be grown starting from $i$ by adding nearest neighbors at each step. 
We call $a^{(i)}$ {\em the $i$-ordering}. 

%\begin{figure}[htbp]%put here figure 5 from picture.pdf
   %\centering
  % \includegraphics{} % requires the graphicx package
  % \caption{\red{An example of an $i$-ordering.}}
   %\label{fig:1}
%\end{figure}
\bigskip

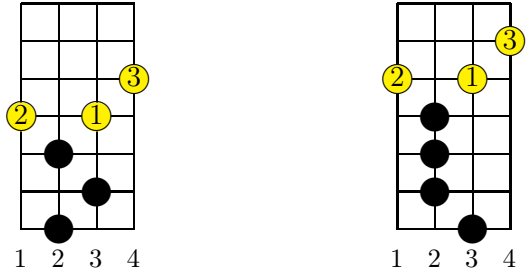
\begin{figure}[htb]
\setlength{\unitlength}{1mm}
\begin{center}
\begin{picture}(85,45)(0,-7)
\multiput(10,5)(5,0){4}{\line(0,1){30}}
\multiput(10,5)(0,5){7}{\line(1,0){15}}
\multiput(60,5)(5,0){4}{\line(0,1){30}}
\multiput(60,5)(0,5){7}{\line(1,0){15}}
\put(15,5){\circle*{4}} 
\put(15,15){\circle*{4}} 
\put(20,10){\circle*{4}} 
\put(10,20){\color{yellow} \circle*{4}} \put(9,19){2} 
\put(10,20){\color{black} \circle{4}}
\put(20,20){\color{yellow} \circle*{4}} \put(19,19){1}   
\put(20,20){\color{black} \circle{4}}
\put(25,25){\color{yellow} \circle*{4}} \put(24,24){3}   
\put(25,25){\color{black} \circle{4}}
\put(70,5){\circle*{4}} 
\put(65,15){\circle*{4}} 
\put(65,10){\circle*{4}} 
\put(65,20){\circle*{4}} 
\put(60,25){\color{yellow} \circle*{4}} \put(59,24){2} 
\put(60,25){\color{black} \circle{4}}
\put(70,25){\color{yellow} \circle*{4}} \put(69,24){1} 
\put(70,25){\color{black} \circle{4}}
\put(75,30){\color{yellow} \circle*{4}}\put(74,29){3} 
\put(75,30){\color{black} \circle{4}}
\put(9,0){\small{1}}
\put(14,0){\small{2}}
\put(19,0){\small{3}}
\put(24,0){\small{4}}
\put(59,0){\small{1}}
\put(64,0){\small{2}}
\put(69,0){\small{3}}
\put(74,0){\small{4}}

\end{picture} 
\end{center} \caption{Once the $i$-ordering is given, from different height profiles having the maximum at the same vertex, using \eqref{3.1} we end up with the same configuration in $S_1$.}
\label{fig:2}
\end{figure}

For the given site $i\in V$ let us write $S^{(i)}=\{y\in S:y_i=0\}$ 
(meaning that the maximum is realized at $i$).  
We put particles according to the corresponding $i$-ordering  $a^{(i)}=(a^{(i)}_{1},\dots,a^{(i)}_{|V|-1})$ and 
look at the resulting configuration
\begin{equation}\label{3.1}
\begin{split}
T_{a^{(i)}_{|V|-1}}\dots T_{a^{(i)}_1}y=:x^{(i)}
\end{split}
\end{equation}
We note that the profile on the  r.h.s. is independent of the choice of $y\in S^{(i)}$ 
and stays bounded with $\min_{j\in V}x^{(i)}_{j}\geq-(|V|-1)$ (see fig. \ref{fig:2}).

Let us put together these configurations and consider the finite subset 
\begin{equation}
\begin{split}
S_{1}=\{ x^{(i)} : i\in |V|\}
\end{split}
\end{equation}
denoting the complement by $S_{2}=S\backslash S_{1}$. 

%\begin{figure}[htbp]%put here figure 2 from picture.pdf
   %\centering
  % \includegraphics{} % requires the graphicx package
   %\caption{{\red Once the $i$-ordering is given, from different height profiles having the maximum at the same vertex, using \eqref{3.1} we end up with the same configuration in $S_1$.}}
   %\label{fig:2}
%\end{figure}
\bigskip

We note the following lemma. 
\begin{lem}
\begin{equation}
\begin{split}\label{killing}
\inf_{x\in S}M^{|V|-1}(x,S_1)\geq \a>0
\end{split}
\end{equation}
where $M^{|V|-1}$ is given by the matrix product. 
\end{lem}
This is clear since any addition of a particle 
has a positive probability and finitely many of those have 
to be considered, leading to the formula $\a=\min_{i\in V}\prod_{j=1}^{|V|-1}p_{a^{(i)}_j}$.  
Next we have the following lemma.

\begin{lem} The equation $\pi M = \pi$
for the invariant distribution has a solution $\pi\in\PP(S)$. 
\end{lem} 

{\bf Proof. } We can say that there is exponential killing
on the infinite part of the space $S_{2}$ and the Markov chain comes
back safely to $S_{1}$. This makes it \textquotedbl{}effectively
finite state\textquotedbl{}. Now, to see this, let 
us introduce the four block-matrices $M_{ij}=(M(x,y))_{x\in S_i, y \in S_j}$, 
introduce the two vectors $\pi_{i}=(\pi(x))_{x\in S_i}$ for $\pi\in \PP(S)$
and rewrite the equation $\pi M = \pi$ for the 
invariant distribution $\pi$ in component form 
\begin{equation}
\begin{split}(\pi_{1},\pi_{2})\left(\begin{array}{cc}
M_{11} & M_{12}\\
M_{21} & M_{22}\end{array}\right)=(\pi_{1},\pi_{2})\end{split}
\end{equation}
 This is equivalent to the form 
  \begin{equation}
\begin{split}\label{2.7} & \pi_{2}=\pi_{1}M_{12}(1_{2}-M_{22})^{-1}\\
 & \pi_{1}(M_{12}(1_{2}-M_{22})^{-1}M_{21}+M_{11})=\pi_{1}\end{split}
\end{equation}
provided that $(1_{2}-M_{22})^{-1}=\sum_{l=0}^{\infty}M_{22}^{l}$ 
exists. But to see the latter use the norm $\Vert M_{2 2} \Vert=\sup_{x \in S_2}\sum_{y\in S_2}M_{2 2}(x,y) 
$ and note that $\Vert M_{2 2}^{|V|-1}\Vert \leq 1-\a $, by \eqref{killing}. 
Hence $(M_{12}(1_{2}-M_{22})^{-1}M_{21}+M_{11})$ is a well-defined
positive matrix on the finite space $S_1$. It is even a stochastic matrix which can be quickly checked 
analytically using the convergence of the geometric sum $\sum_{l=0}^{\infty}M_{22}^{l}$. 
So the matrix has a Perron-Frobenius eigenvector to the eigenvalue $1$, which we call $\pi_{1}^{*}$
(up to a positive multiple). 
This is (up to this multiple) the invariant distribution restricted to $S_{1}$.
From this we get the invariant distribution $\pi^*_2$ on the infinite part of
the system by looking at the first equation of \eqref{2.7} and normalizing. $\Cox$
\bigskip

{\bf Remark:} If $S_1=\{x\}$ is a single point then define the return time 
$\t_x=\inf\{t\geq 1: W(t)=x\}$ where $W(t)$ is a random walk started 
at $x$. For a state $y\neq x$ we have that the non-normalized 
distribution at $y$ is given by the expected number of visits 
from $x$ to $y$ before returning to $y$, i.e.   
$[M_{12}(1_{2}-M_{22})^{-1}]_{x,y}=\E^x \sum_{t=1}^\infty 
1_{W(t)=y}1_{t <\t_x}$. Normalization of the distribution 
then implies that $1=\pi(x)+\pi(x)(\E^x \t_x -1)$  and so $\pi(x)=\frac{1}{\E^x \t_x}$ and $\pi(y)=\frac{1}{\E^x \t_x}\E^x \sum_{t=1}^\infty 
1_{W(t)=y}1_{t <\t_x}$ for $y\neq x$.

\begin{lem} The Markov chain is {\em uniformly communicating} to $S_1$ by which 
we  mean that there exists 
an $\a'>0$ and a time $s= 3(|V|-1)$, called the {\em communication time}, such that 
\begin{equation}
\begin{split}\label{communism}
\inf_{x\in S, x' \in S_1}M^{s}(x,x')\geq \a'>0
\end{split}
\end{equation} 
\end{lem}

{\bf Proof. } The proof follows by noting that we can first: get 
into $S_1$, second: go from there into a state
which has a prescribed maximum (possibly outside $S_1$), 
and third: go from that 
state into the corresponding state in $S_1$. In formulas it reads like this: 
Consider a starting configuration $y\in S^{(i)}$. Then, with 
the above construction we have 
\begin{equation}
\begin{split}
T_{a^{(i)}_{|V|-1}}\dots T_{a^{(i)}_1}y=x^{(i)}\in S_1
\end{split}
\end{equation}
We note that $(T_{j})^{|V|-1}x^{(i)}\in S^{(j)}$ since sufficiently many 
particle droppings at $j$ are shifting the maximum to the point $j$. 
From that we get again by the first step that 
\begin{equation}
\begin{split}
x^{(j)}=T_{a^{(j)}_{|V|-1}}\dots T_{a^{(j)}_1}(T_{j})^{|V|-1}x^{(i)}
=T_{a^{(j)}_{|V|-1}}\dots T_{a^{(j)}_1}(T_{j})^{|V|-1}T_{a^{(i)}_{|V|-1}}\dots T_{a^{(i)}_1}y
\end{split}
\end{equation}
The proof is complete since $j\in V$ was arbitrary.
$\Cox$
\bigskip

{\bf Remark. } From the above definition of a communication set 
$S_1$ follows trivially that any subset is also a communication set 
since the $\inf$ has to be taken over less terms.  
While from a theoretical point of view it would be 
therefore sufficient to consider a single point $x_0\in S_1$ in our example 
for our chain, returns are easiest understood when we talk about 
our definition of $S_1$.  {The remark will be clear after dealing with particle droppings according to a Markov chain (see fig. \ref{fig:4} and fig. \ref{fig:5}).}

We have from this the convergence to the invariant distribution 
in total variation: 
\begin{lem} 
\begin{equation}
\begin{split}
\Vert M^{s}(x, \cdot) -\pi \Vert_{\text{TV}}\leq 1-(\a')^2 |V|
\end{split}
\end{equation}
\end{lem}

{\bf Proof. }
Call $X_t$ the chain starting at $x$ and $Y_t$ the one starting with initial distribution $\pi$. Moreover call $\tau$ the random time of their first meeting in 
the product coupling. After they meet for the first time they stay together.
The coupling inequality gives:
\begin{multline}\label{1} || M^s(x,\cdot)-\pi||_{TV}\leq\mathbb{P}_c(\tau>s)\leq 1-\sum_{y\in S}\mathbb{P}_x(X_s=y)\mathbb{P}_{\pi}(Y_s=y)\\
=1-\sum_{y\in S}\mathbb{P}_x(X_s=y)\pi(y)\leq 1-\sum_{y\in S_1}\underbrace{\mathbb{P}_x(X_s=y)}_{\geq\alpha'}\underbrace{\pi_1(y)}_{\geq\alpha'}\leq 1-(\a')^2 |V|.
\end{multline} 
 $\Cox$
\bigskip

From the Lemma follows the convergence by standard arguments, 
extending the Lemma to 
$\Vert M^{s k}(x, \cdot) -\pi \Vert_{\text{TV}}\leq (1-(\a')^2 |V|)^k$ for integer $k$ 
and using that the total variation distance is decreasing in the time $t$.

\subsection{Concentration properties of path functionals}

Define, for $t'>t$, the coupling matrix 
\begin{equation}
\begin{split} &
D_{t,t'} := \sup_{x,x'} \P_c( X(t') \neq X'(t') | X(t)=x, X'(t)=x')
\end{split}
\end{equation}
where $\P_c$ is the product coupling mentioned above. 
We have for times which differ by the communication time $s$ that 
\begin{equation}
\begin{split} &
D_{t,t+s}\leq 1-(\a')^2 |V|
\end{split}
\end{equation}
and this implies for general times 
\begin{equation}
\begin{split} &
D_{t,t'} \leq (1-(\a')^2 |V|)^{\lfloor\frac{ t'-t }{s}\rfloor}
\end{split}
\end{equation}

\begin{lem}\label{concentration}
Let $g: S^n \rightarrow \R$ be a bounded measurable function. 
Then we have a Gaussian concentration bound of the form 
\begin{equation}
\begin{split} &
\P(|g -\E g | > y) \leq 2\exp\Bigl( 
- \frac{2 y^2 }{\Vert D \d g \Vert^2 }
\Bigr) 
\end{split}
\end{equation}
where 
\begin{equation}
\begin{split} &
(\d g)_u = \d_u (g) = \max_{x_u, x'_u} g(x_1, \dots, x_u, \dots ) -  g(x_1, \dots, x'_u, \dots )
\end{split}
\end{equation}
is the variation at the time $u$.  
\end{lem}

{\bf Proof. } In the following we give only the key steps in the proof of lemma \eqref{concentration}. We refer to \cite{CCKR06} for details where the same proof in the context of models with finite state space was given. This is not a problem 
here since our observable is bounded, and, most importantly the coupling matrix $D$ satisfies the nice bounds given above, 
in spite of our state space being unbounded, due to uniform coupling speed. 
\\
Using the standard decomposition into 
Martingale differences and the Markov property, we can write 
\begin{equation}
\begin{split} &
g -\E g =\sum_{i=1}^n W_i(x_{1},\dots, x_i) 
\end{split}
\end{equation}
with  
\begin{equation}\label{m1}
\begin{split} &
W_i(x_{1},\dots, x_i) = 
\E(g |x_{1},\dots, x_i)- \E(g |x_{1},\dots, x_{i-1})\cr
%&=\int \P( d \tilde x_{i+1}, \dots ,d \tilde x_{n}| x_i) g(x_1, \dots, x_i, \tilde x_{i+1}, \dots ,\tilde x_{n})\cr
%&-\int \P( d \tilde x_{i}, \dots ,d \tilde x_{n}| x_{i-1}) g(x_1, \dots, x_{i-1}, \tilde x_{i}, \dots ,\tilde x_{n})\cr
%&=\int \P( d \tilde x_{i+1}, \dots ,d \tilde x_{n}| x_i) g(x_1, \dots, x_i, \tilde x_{i+1}, \dots ,\tilde x_{n})\cr
%&-\int \P( d \tilde x_{i}| x_{i-1}) \int \P( d \tilde x_{i+1}, \dots ,d \tilde x_{n}| x_i) g(x_1, \dots, x_{i-1}, \tilde x_{i}, \dots ,\tilde x_{n})\cr
&\leq \sup_{\bar x_i\in S }\int \P( d \tilde x_{i+1}, \dots ,d \tilde x_{n}| \bar x_i) g(x_1, \dots, \bar x_i, \tilde x_{i+1}, \dots ,\tilde x_{n})\cr
&- \inf_{\bar y_i\in S }
\int \P( d \tilde x_{i+1}, \dots ,d \tilde x_{n}| \bar y_i) g(x_1, \dots, \bar y_i, \tilde x_{i+1}, \dots ,\tilde x_{n})\cr
&=:X_i(x)-Y_i(x)\cr
\end{split}
\end{equation}
Notice that the $\inf$ and $\sup$ appearing in the previous formula are well defined since $g$ is bounded.
Then, we use a simple telescoping identity to rewrite $g$ as a sum of discrete gradients 
 \begin{equation}
\begin{split} &g(x_1, \dots, \bar x_i, \tilde x^{(1)}_{i+1}, \dots ,\tilde x^{(1)}_{n})-g(x_1, \dots, \bar y_i, \tilde x^{(2)}_{i+1}, \dots ,\tilde x^{(2)}_{n})\cr
%&=g(x_1, \dots, \bar x_i, \tilde x^{(1)}_{i+1}, \dots ,\tilde x^{(1)}_{n})-g(x_1, \dots, \bar y_i, \tilde x^{(1)}_{i+1}, \dots ,\tilde x^{(1)}_{n})\cr
%&+g(x_1, \dots, \bar y_i, \tilde x^{(1)}_{i+1}, \dots ,\tilde x^{(1)}_{n})-g(x_1, \dots, \bar y_i, \tilde x^{(2)}_{i+1}, \tilde x^{(1)}_{i+2}\dots ,\tilde x^{(1)}_{n})\cr
%&+\ldots\cr
%&+g(x_1, \dots, \bar y_i, \tilde x^{(2)}_{i+1}, \dots ,\tilde x^{(2)}_{n-1},\tilde x^{(1)}_{n})-g(x_1, \dots, \bar y_i, \tilde x^{(2)}_{i+1}, \tilde x^{(2)}_{i+2}\dots ,\tilde x^{(2)}_{n})\cr
&=\sum_{j=0}^{n-i}\nabla_{i,i+j}^{12}g.
\end{split}
\end{equation}
where $\nabla_{i,i+j}^{12}g$ is the difference between $g$'s evaluated at two points that are the same except for the $(i+j)-th$ place. We define:
\begin{equation}
(\d g)_u = \d_u (g) = \max_{x_u, x'_u} g(x_1, \dots, x_u, \dots ) -  g(x_1, \dots, x'_u, \dots ),
\end{equation}
and by construction we have
\begin{equation}\label{m2}
\nabla_{i,i+j}^{12}g\leq  \d_{i+j} (g)1_{\tilde x^{(1)}_{i+j}\neq\tilde x^{(2)}_{i+j} }.
\end{equation}
Then using \eqref{m1} and \eqref{m2}, it follows that
 \begin{equation}\label{m3}
\begin{split} X_i(x)-Y_i(x)&=\sup_{\bar x_i, \bar y_i\in S}\Bigl\{\int \P( d \tilde x_{i+1}, \dots ,d \tilde x_{n}| \bar x_i) g(x_1, \dots, \bar x_i, \tilde x_{i+1}, \dots ,\tilde x_{n})\cr
&\P( d \tilde x_{i+1}, \dots ,d \tilde x_{n}| \bar y_i) g(x_1, \dots, \bar y_i, \tilde x_{i+1}, \dots ,\tilde x_{n})\Bigr\}\cr
&=\sup_{\bar x_i, \bar y_i\in S}\Bigl\{\int \P_c(d \tilde x^{(1)}_{\geq i+1}, d \tilde x^{(2)}_{\geq i+1}|\tilde x^{(1)}_i=\bar x_i, \tilde x^{(2)}_i=\bar y_i)[g(x_1, \dots, \bar x_i, \tilde x^{(1)}_{i+1}, \dots ,\tilde x^{(1)}_{n})\cr
&-g(x_1, \dots, \bar x_i, \tilde x^{(2)}_{i+1}, \dots ,\tilde x^{(2)}_{n})]\Bigr\}\cr
&\leq \sup_{\bar x_i, \bar y_i\in S}\sum_{j=0}^{n-i} \d_{i+j} (g)\P_c\left(\tilde x^{(1)}_{i+j}\neq\tilde x^{(2)}_{i+j} |\tilde x^{(1)}_i=\bar x_i,\tilde x^{(2)}_i=\bar y_i\right)\cr
&\leq \sum_{j=0}^{n-i}  D_{i,i+j}\d _{i+j}g=:(D\d g)_i.\cr
\end{split}
\end{equation}
The last ingredient is the following lemma from \cite{DeLu01}.
\begin{lem} Suppose $\mathcal {F}$ is a $\sigma$-field and $Z_1, Z_2, W$ are random variables  such that
\begin{enumerate}
\item $Z_1\leq W \leq Z_2$;
\item $\mathbb{E}(W|\mathcal {F})=0$;
\item $Z_1$ and $Z_2$ are $\mathcal{F}$-measurable.
\end{enumerate}
Then, for all $\lambda\in\mathbb{R}$, we have the inequality
\begin{equation}
\mathbb{E}\left[\exp(\l W)|\mathcal{F}\right]\leq \exp\left[\frac{\l^2(Z_2-Z_1)^2}{8}\right].
\end{equation}
\end{lem} 
This lemma, in the present situation, works putting $W=W_i$,  $Z_1=X_i-\mathbb{E}[g|\mathcal{F}_{i-1}]$, $Z_2~=~Y_i~-~\mathbb{E}[g|\mathcal{F}_{i-1}]$ and $\mathcal F=\mathcal F_{i-1}$. Since, from \eqref{m3} and \eqref{m1} we have
\begin{equation}
W_i\leq Y_i-X_i\leq (D\d g)_i,
\end{equation}
we obtain
\begin{equation}\label{m4}
\mathbb{E}\left[\exp(\l W_i)|\mathcal{F}_{i-1}\right]\leq \exp\left[\frac{\l^2(D\d g)_i}{8}\right].
\end{equation}
By the exponential Chebyshev inequality and iterating \eqref{m4} by successive conditional expectations with respect to $\mathcal{F}_n$ we compute
\begin{equation}
\begin{split}\mathbb{P}(g-\mathbb{E}g\geq y)
%&\leq \exp[-\l y]\mathbb{E}\left[\exp\left(\l\sum_{i=1}^{n}W_i\right)\right]\cr
%&\leq \exp[-\l y]\mathbb{E}\left\{\mathbb{E}\left[\exp(\l W_n)|\mathcal{F}_{n-1}\right]\exp\left(\l\sum_{i=1}^{n-1}W_i\right)\right\}\cr
&\leq\exp[-\l y]\exp\left[\frac{\l ^2}{8}||D\d g||^2\right].\cr
\end{split}
\end{equation}
We choose the optimal $\l=4y/||D\d g||^2$ to obtain
\begin{equation}
\mathbb{P}(g-\mathbb{E}g\geq y)\leq\exp\left[ -\frac{2y^2}{||D\d g||^2}\right].
\end{equation}
The previous line of reasoning applies to $-g$ and $-W$, proving \eqref{concentration}.   
$\Cox$\bigskip\bigskip

\subsection{The total height as an additive path functional}

Let us come back now to our main application and consider the maximum of 
the total height of the original process, started from the flat initial configuration at zero, given by
\begin{equation}m_V(t)=\max_{j\in V}h_{j}(t)\end{equation}
The main idea is to write a formula as an additive functional of the Markov chain along the path: 
\begin{equation}
\begin{split} &
m_V(t)=t-\sum_{u=1}^{t-1}1_{\max_j h_j(u+1)=\max_j h_j(u)}
\end{split}
\end{equation}
We will rewrite the functions under the sum in terms of the $x$-process instead of the 
original one using the following lemma using the following numbers. 
\begin{lem} 
\begin{equation}
\begin{split} &
\max_j h_j(u+1)=\max_j h_j(u) \Leftrightarrow \#\{j \in V | x_j(u)\neq x_j(u+1)  \}=1
\end{split}
\end{equation}
\end{lem}
{\bf Proof. } To see that the l.h.s. implies the r.h.s. 
note that under the assumption $\max_j h_j(u+1)=\max_j h_j(u)$ we have 
$\#\{j \in V | x_j(u)\neq x_j(u+1)  \}=\#\{j \in V | h_j(u)\neq h_j(u+1)\}=1$

To see that the r.h.s. implies the l.h.s., let us suppose that $\max_j h_j(u+1)=\max_j h_j(u)+~1$,
and derive a contradiction. But indeed in that case we would have 
$\#\{j \in V | x_j(u)\neq x_j(u+1)  \}=\#\{j \in V | h_j(u)\neq h_j(u+1)-1\}=|V|-1$ 
which is different from $1$ if $|V|>2$.  
$\Cox$
\bigskip

So we have 
\begin{equation}
\begin{split} &
m_V(t)=t-\sum_{u=1}^{t-1}1_{ \#\{j \in V | x_j(u)\neq x_j(u+1)  \}=1}
\end{split}
\end{equation}
In our case we have $\d_v (\sum_{u=1}^{t-1}1_{ \#\{j \in V | x_j(u)\neq x_j(u+1)  \}=1}) \leq 2$
giving us 
\begin{equation}
\begin{split} 
\P(m_V(t) -\E m_V(t) > y) &\leq \exp\Bigl( 
- \frac{y^2 }{2  \sum_{u=1}^t (\sum_{u':u< u' <t }D_{u,u'})^2  } \Bigr)
\leq \exp\Bigl( 
- \frac{y^2 }{2 t} \frac{|V|^2(\a')^4}{s^2} 
\Bigr) 
\end{split}
\end{equation}
and the same bound for $\P(m_V(t) -\E m_V(t) <- y) $.
Both bounds hold if we take for $\P=\P_{\pi}$ the chain in equilibrium 
or with a given initial condition, say $h=0$. Denote this chain by $\P_0$.  
Since our original interest was in the latter one we need to note 
the closeness of the two expected values  which follows again by using the 
uniform bound on the coupling to compare the two distributions in the second inequality of 
\begin{equation}
\begin{split} 
&\Bigl |\E_{\pi} m_V(t)- \E_0 m_V(t)\Bigr|\cr 
&\leq \sum_{u=1}^{t-1}\Bigl |
 \E_0 \Bigl(1_{  \#\{j \in V | x_j(u)\neq x_j(u+1)  \}=1} \Bigr)
- \E_{\pi} \Bigl(1_{ \#\{j \in V | x_j(u)\neq x_j(u+1)  \}=1}\Bigr)\Bigr|\cr
 &\leq \sum_{u=1}^{t-1} D_{0,u}\leq \frac{s}{|V|(\a')^2} 
\end{split}
\end{equation}
In particular we get the strong law of large numbers 
\begin{equation}
\begin{split} 
&\lim_{t\uparrow \infty}\frac{m_V(t)}{t}= 1-\sum_{x,y}\pi(x)M(x,y)1_{ \#\{j \in V | x_j\neq y_j  \}=1}
\end{split}
\end{equation}
\bigskip 

{ This is a particular example of an Ergodic Theorem for path observables which enjoy the concentration property.}

\section{Markov chain particle droppings - the proof }

We consider the mapping from $h(t)$ to $x(t)$ as above and remark 
that $(x(t),v(t))$ is a Markov chain again. 
Warning: It is not to be expected that the marginal process $x(t)$ 
is a Markov chain  (of  memory depth $1$) now. It will be a chain with a depth of memory $2$ since 
the position of $v(t)$ can be reconstructed looking by $(x(t),x(t-1))$.

It is useful to make explicit the graph $(V,E)$ with undirected edges $E$ defining the piling-up rule, 
and the graph $(V,E_A)$ with directed edges $E_A=\{ (i,j)\in V\times V: A(i,j)>0  \}$. 
The following considerations depend on $A$ only through $E_A$. 

%\begin{figure}[htbp]%put here figure 2 from picture.pdf
   %\centering
  % \includegraphics{} % requires the graphicx package
   %\caption{{\red Not all edges in the graph $(V,E)$ (on the left) correspond to two 
  % directed edges in the graph $(V,E_A)$. 
  %The graphs $(V,E)$ (on the left) and $(V,E_A)$ (on the right) may be different. 
  %This may forbid us to put particles neighboring each other in one step.}}
   %\label{fig:3}
%\end{figure}

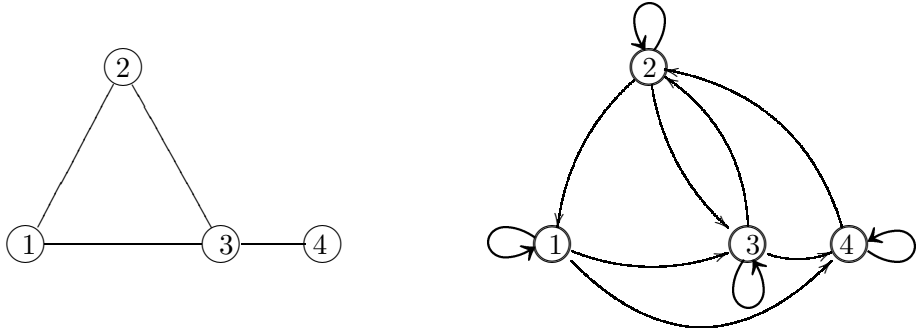
\begin{figure}[htb]
\setlength{\unitlength}{1mm}
\begin{center}
\begin{picture}(125,45)
\put(7.5,35){\xymatrix{
           & 2 \ar@{-}[ldd]  \ar@{-}[rdd] &   &  \\
           &                    &   &  \\  
1  \ar@{-}[rr]   &                    & ~3 \ar@{-}[r] & 4 \\
}}

\multiput(22,35)(69.9,0){2}{\circle{5}}
\multiput(9,11.5)(70,0){2}{\circle{5}}
\multiput(35,11.5)(70,0){2}{\circle{5}}
\multiput(48.5,11.5)(70,0){2}{\circle{5}}

\put(91.9,35){\circle{4.5}}
\put(79,11.5){\circle{4.5}}
\put(105,11.5){\circle{4.5}}
\put(118.5,11.5){\circle{4.5}}

\put(77.5,35){\xymatrix{
           & 2  \ar@/^-0.3cm/[rdd] \ar@/^-0.3cm/[ldd] &   &  \\
           &                    &   &  \\  
1 \ar@/^-1.1cm/[rrr] \ar@/^-0.3cm/[rr]   &  & ~3 \ar@/^-0.2cm/[r] \ar@/^-0.4cm/[luu] & 4 \ar@/^-0.6cm/[lluu]\\
}}
\put(70.0,9.9){\huge{$\lefttorightarrow$}}
\put(102.2,4.3){\huge{$\downtouparrow$}}
\put(89.1,38.6){\huge{$\uptodownarrow$}}
\put(120.5,9.4){\huge{$\righttoleftarrow$}}
\end{picture} 
\end{center} \caption{Not all edges in the graph $(V,E)$ (on the left) correspond to two 
   directed edges in the graph $(V,E_A)$. 
  %The graphs $(V,E)$ (on the left) and $(V,E_A)$ (on the right) may be different. 
  This may forbid us to put particles neighboring each other in one step.}
\label{fig:3}
\end{figure}
\bigskip

The first step is to extend the definition of $S_1$ to the present setup. 
A problem might be that the driving Markov chain forbids us to put 
balls neighboring each other {\em in one step} in the sense of the graph which defines our piling-up rules (see fig. \ref{fig:3}). 
What we need is to make sure that we can define a finite set $\bar S_1$ to which the joint chain 
communicates uniformly. While in the independent case we could just take the configurations
which were grown from nearest neighbor to nearest neighbor (along an $i$-ordering), 
here we have to add to it connecting strings of allowed transitions in between. A slight 
discomfort is that the maximum might change in a complicated way during this process of particle 
additions following this string.  
However, this is not really important. What is important is that a resulting configuration 
will only depend on the maximum of the initial configuration and otherwise be independent on its form. 
Now, we can ensure the latter by adding sufficiently many particles at the maximizing site initially.  
By lazyness it is a chain of allowed transitions
and it ensures that all influence of the configuration at any other site will be lost. 
This is formulated in the following Lemmata. The first Lemma is treating 
a situation where the driving Markov chain drops a particle at the same position as the maximum. 
The next Lemma shows how the situation where the driving Markov chain is in 
a different position than the maximum can be reduced to the first. 

\begin{lem} \label{4.1}Suppose that $x\in S^{(i)}$. Then there exists a finite integer  
$s(i)$ and a sequence $(i=i_1,i_2,\dots, i_{s(i)})\in V^{s(i)}$ such that $(i_j,i_{j+1})\in E_A$ 
is an allowed transition and the configuration
\begin{equation}
\begin{split}
x^{(i)}:=T_{i_{s(i)}} \dots T_{i_1}(T_i)^{s(i)} x 
\end{split}
\end{equation}
is independent of the choice of the initial configuration in $S^{(i)}$ and 
has a bounded depth $\min_{j\in V}x^{(i)}_j\geq - 2 s(i)$. 
\end{lem}

\begin{figure}[htb]
\setlength{\unitlength}{1mm}
\begin{center}
\begin{picture}(120,85)(0,-7)

\put(20,80){\Large{(a)}}
\put(89,80){\Large{(b)}}
\put(28,4){\Large{(c)}}

\put(40,5){\xymatrix{
1 \ar@/^0.4cm/[rr] & 2 \ar@/^0.4cm/[rr] & 3 \ar@/^/[l] & 4 \ar@/^-0.9cm/[lll] \\
}}
\multiput(38.9,-1.8)(12.5,0){4}{\huge{$\downtouparrow$}}
\multiput(41.8,5)(12.5,0){4}{\circle{4}}
\multiput(41.8,5)(12.5,0){4}{\circle{3.5}}
\multiput(0,25)(5,0){4}{\line(0,1){50}}
\multiput(0,25)(0,5){11}{\line(1,0){15}}

\multiput(30,25)(5,0){4}{\line(0,1){15}} \multiput(30,55)(5,0){4}{\line(0,1){20}}
\multiput(30,25)(0,5){4}{\line(1,0){15}} \multiput(30,55)(0,5){5}{\line(1,0){15}}
%\multiput(29.5,26)(5,0){4}{\vdots}

\multiput(70,25)(5,0){4}{\line(0,1){50}}
\multiput(70,25)(0,5){11}{\line(1,0){15}}

\multiput(100,25)(5,0){4}{\line(0,1){15}} \multiput(100,55)(5,0){4}{\line(0,1){20}}
\multiput(100,25)(0,5){4}{\line(1,0){15}} \multiput(100,55)(0,5){5}{\line(1,0){15}}

\put(5,25){\circle*{4}} 
\put(0,30){\circle*{4}} 
\put(10,30){\circle*{4}} 
\put(5,35){\circle*{4}} 
\put(15,25){\circle*{4}} 
\put(35,25){\circle*{4}} 
\put(30,30){\circle*{4}} 
\put(40,30){\circle*{4}} 
\put(35,35){\circle*{4}} 
\put(45,25){\circle*{4}} 

%\put(0,20){$1$}
\put(0,40){\color{yellow}  \circle*{4}}  
\put(0,40){\color{black}  \circle{4}} \put(-1,39){2} 
\put(10,40){\color{yellow} \circle*{4}}
\put(10,40){\color{black}  \circle{4}} \put(9,39){3}   
\put(15,35){\color{yellow} \circle*{4}} 
\put(15,35){\color{black}  \circle{4}} \put(14,34){1} 
\put(5,45){\color{yellow}  \circle*{4}}
\put(5,45){\color{black}  \circle{4}} \put(4,44){4} 
\put(15,45){\color{yellow} \circle*{4}}
\put(15,45){\color{black}  \circle{4}} \put(14,44){5} 

\put(35,40){\circle{4}} 
\put(35,40){\color{cyan} \circle*{4}} \put(35,40){\color{black} \circle*{2}}  
\put(34.5,46){\vdots}
\put(35,55){\circle{4}} 
\put(35,55){\color{cyan} \circle*{4}} \put(35,55){\color{black} \circle*{2}} 

\put(42,46.6){\Large{$s(i)$}}
\put(33,46.6){$\left. \begin{array}{r} \\ \\ \\ \end{array} \right\}$ } 

\put(30,60){\color{yellow} \circle*{4}} 
\put(30,60){\color{black}  \circle{4}} \put(29,59){2} 
\put(40,60){\color{yellow} \circle*{4}} 
\put(40,60){\color{black}  \circle{4}} \put(39,59){3} 
\put(45,65){\color{yellow} \circle*{4}} 
\put(45,65){\color{black}  \circle{4}} \put(44,64){5} 
\put(35,65){\color{yellow} \circle*{4}}
\put(35,65){\color{black}  \circle{4}} \put(34,64){4} 
\put(45,35){\color{yellow} \circle*{4}}
\put(45,35){\color{black}  \circle{4}} \put(44,34){1} 

\put(70,25){\circle*{4}} \put(80,25){\circle*{4}} 
\put(80,30){\circle*{4}} \put(75,35){\circle*{4}}
\put(85,35){\circle*{4}}

\put(70,40){\color{yellow} \circle*{4}} 
\put(70,40){\color{black} \circle{4}} \put(69,39){2} 
\put(85,40){\color{yellow} \circle*{4}}
\put(85,40){\color{black} \circle{4}} \put(84,39){1}  
\put(80,45){\color{yellow} \circle*{4}} 
\put(80,45){\color{black} \circle{4}} \put(79,44){3} 
\put(75,50){\color{yellow} \circle*{4}}
\put(75,50){\color{black} \circle{4}} \put(74,49){4} 
\put(85,50){\color{yellow} \circle*{4}}
\put(85,50){\color{black} \circle{4}} \put(84,49){5} 
\put(100,25){\circle*{4}} \put(110,25){\circle*{4}} 
\put(110,30){\circle*{4}} \put(105,35){\circle*{4}}
\put(115,35){\circle*{4}}

\put(105,40){\circle{4}}
\put(105,40){\color{cyan} \circle*{4}} \put(105,40){\color{black} \circle*{2}} 
\put(104.5,46){\vdots}
\put(105,55){\circle{4}} 
\put(105,55){\color{cyan} \circle*{4}} \put(105,55){\color{black} \circle*{2}} 

%\put(108,48){$\backslash$} \put(108,44.5){$/$} 
\put(112,46.6){\Large{$s(i)$}}
\put(103,46.6){$\left. \begin{array}{r} \\ \\ \\ \end{array} \right\}$ } 

\put(100,60){\color{yellow} \circle*{4}} 
\put(100,60){\color{black} \circle{4}} \put(99,59){2} 
\put(115,65){\color{yellow} \circle*{4}} 
\put(115,65){\color{black} \circle{4}} \put(114,64){5} 
\put(110,60){\color{yellow} \circle*{4}} 
\put(110,60){\color{black} \circle{4}} \put(109,59){3} 
\put(105,65){\color{yellow} \circle*{4}}
\put(105,65){\color{black} \circle{4}} \put(104,64){4} 
\put(115,40){\color{yellow} \circle*{4}}
\put(115,40){\color{black} \circle{4}} \put(114,39){1} 
\put(-1,20){\small{1}}
\put(4,20){\small{2}}
\put(9,20){\small{3}}
\put(14,20){\small{4}}
\put(29,20){\small{1}}
\put(34,20){\small{2}}
\put(39,20){\small{3}}
\put(44,20){\small{4}}

\put(69,20){\small{1}}
\put(74,20){\small{2}}
\put(79,20){\small{3}}
\put(84,20){\small{4}}

\put(99,20){\small{1}}
\put(104,20){\small{2}}
\put(109,20){\small{3}}
\put(114,20){\small{4}}

\end{picture} 
\end{center} \caption{{ An example of the procedure described in the proof of Lemma \ref{4.1} on the graph with edges $\{\{1,2\},\{2,3\},\{3,4\}\}$.  Dropping particles according to the concatenated strings 
   along the $i=2$-ordering $(1,3,4)$  does not suffice to obtain the same configurations in $\bar S_1$ (see the left pictures of $(a)$ and $(b)$).  If we previously add $s(i)$ balls to the top of the height profiles the configurations are the same (see the right pictures of $(a)$ and $(b)$).}}
\label{fig:4}
\end{figure}
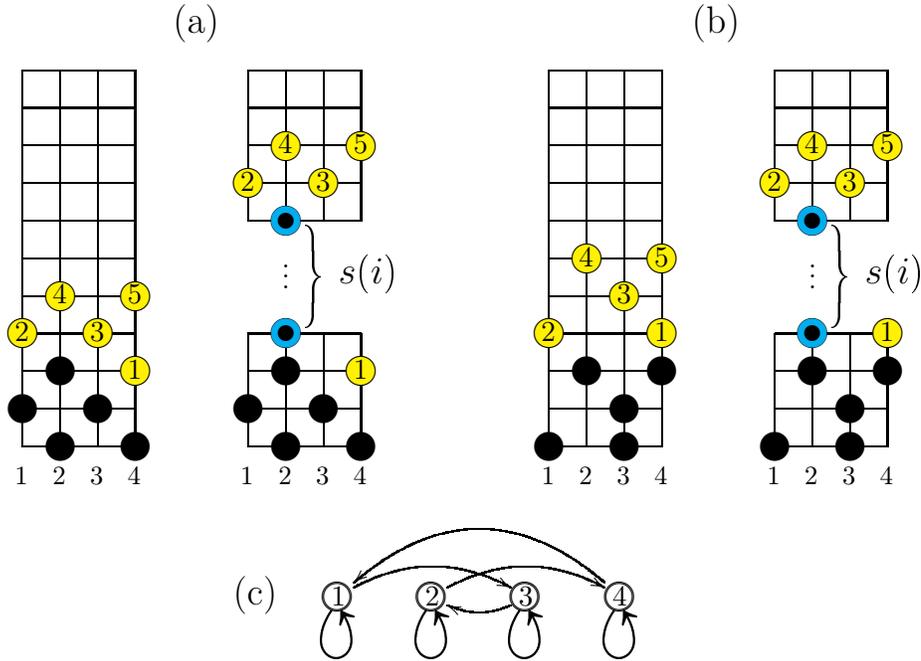

{\bf Proof. }
We choose for each vertex $i\in V$ an $i$-ordering 
$a^{(i)}=(a^{(i)}_{1},\dots,a^{(i)}_{|V|-1})$
of the sites in the set $V\backslash \{i\}$, which was 
defined above. 
We need to connect each of the occurring pairs of neighboring vertices $v=a^{(i)}_{j}, 
w=a^{(i)}_{j+1}$
with a chain of allowed transitions $(v_1=v,v_2 ,\dots, v_{s(v,w)}=w)$ where $s(v,w)$ is the shortest 
length of an oriented path in $E_A$. 
In particular every vertex in the string is visited only once.  Let us denote 
the string from $v$ to $w$ which we obtain by the above 
by dropping the $w$ from it by  $c(v,w)=(v_1,\dots, v_{s(v,w)-1})$. 
Then we concatenate the strings along the $i$-ordering and define 
$$(i_1,\dots, i_{s(i)}):=(c(i,a^{(i)}_{1}), c(a^{(i)}_{1},a^{(i)}_{2}),c(a^{(i)}_{2},a^{(i)}_{3})),\dots, c(a^{(i)}_{|V|-2},
a^{(i)}_{|V|-1}),a^{(i)}_{|V|-1})$$
This string has the property that it contains the $i$-ordering as a substring and 
therefore erases the influence of an initial configuration $y\in S^{(i)}$ when applied to it, when 
the difference of the maximum at $i$ and the configuration 
at any other site was bigger than any possible number of occurrences 
of a site $j$ in  $(i_1,\dots, i_{s(i)})$ (see fig\ref{fig:4}).  
$\Cox$
\bigskip 
\bigskip
%\begin{figure}[htbp]%put here figure 3 from picture.pdf
   %\centering
  % \includegraphics{} % requires the graphicx package
   %\caption{{\red An example of the procedure described in the proof of Lemma \ref{4.1}.  Dropping particles according to the concatenated strings 
   %along the $i$-ordering $(2,1,3,4)$  does not suffice to obtain the same configurations in $\bar S_1$ (see the left pictures of $(a)$ and $(b)$).  If we previously add $s(i)$ balls to the top of the height profiles the configurations are the same (see the right pictures of $(a)$ and $(b)$).}}
   %\label{fig:4}
%\end{figure}

\bigskip

\begin{lem} Suppose that $x\in S^{(i)}$ and $v\in V$. 
Then there exists a finite integer  
$\s(i)$ and a sequence $(i_1=v,i_2,\dots, i_{\s(i)}=i)\in V^{\s(i)}$ such that $(i_j,i_{j+1})\in E_A$ 
is an allowed transition and
\begin{equation}
\begin{split}
T_{i_{\s(i)}} \dots T_{i_1}x \in S^{(i)}
\end{split}
\end{equation}
\end{lem}

The Lemma says we can go from any initial position of the driving Markov chain 
and a height profile with maximum in $i$ to a position with maximum again in $i$ 
and driving Markov chain also in $i$, just as the first lemma assumed.\medskip

{\bf Proof. } First drop $s(v,i)$ particles according to $c(v,i)$. Then drop $s(v,i)$ particles at $i$ 
to be sure that the maximum will be again at $i$. 
This proves the lemma with $\s(i)=2 \max_v s(v,i)$. 
$\Cox$
\bigskip

In analogy to the independent case we put together these configurations and consider the finite subset 
(which this time however will live in the product space $\bar S=S\times V$) and define 
\begin{equation}
\begin{split}
\bar S_{1}=\{ (x^{(i)},i_{s(i)}) : i\in |V|\}
\end{split}
\end{equation}
denoting the complement by $\bar S_{2}=\bar S\backslash \bar S_{1}$. 

\begin{corollar}\label{4.3} The Markov chain $M$ is uniformly communicating to $\bar S_{1}$ 
with a finite communication time $\bar s$.
%$=3 \t + \max_{v\neq w} s(v,w)$ where $\t:=\max_{i\in V}2 s^{(i)}$,  
%and $\bar \a'>0$. 
\end{corollar}

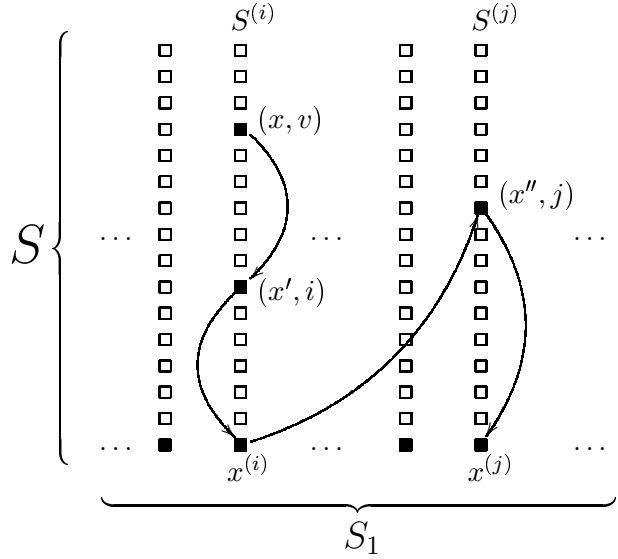
\begin{figure}[htb]
\setlength{\unitlength}{1mm}
\begin{center}
\begin{picture}(80,80)(0,-7)

\put(30,55){\xymatrix{
\ar@/^0.6cm/[dd] &   &  &  &  \\
                 &   &  & \ar@/^0.6cm/[ddd] &                  \\
\ar@/^-0.6cm/[dd]&   &  &  &                  \\
                     &  &  &  &                  \\
\ar@/^-0.6cm/[rrruuu]&  &  &  &                  \\
}}
\multiput(20,12)(10,0){2}{\multiput(0,0)(0,3.5){16}{\line(1,0){1.5} \line(0,1){1.5} }}
\multiput(20,13.5)(10,0){2}{\multiput(0,0)(0,3.5){16}{\line(1,0){1.5}}}
\multiput(20,13.5)(10,0){2}{\multiput(0,0)(0,3.5){16}{\line(0,-1){1.5}}}

\multiput(52,12)(10,0){2}{\multiput(0,0)(0,3.5){16}{\line(1,0){1.5} \line(0,1){1.5} }}
\multiput(52,13.5)(10,0){2}{\multiput(0,0)(0,3.5){16}{\line(1,0){1.5}}}
\multiput(52,13.5)(10,0){2}{\multiput(0,0)(0,3.5){16}{\line(0,-1){1.5}}}

\put(12,40){\dots}
\put(40,40){\dots}
\put(75,40){\dots}

\put(12,12){\dots}
\put(40,12){\dots}
\put(75,12){\dots}
\put(0,37){\huge{$S$}}
\put(5,38){$\left\{ \begin{array}{r} \\ \\ \\ \\ \\ \\ \\ \\ \\ \\ \\ \\  
\end{array} \right.$ } 

\put(12,6){$\underbrace{~~~~~~~~~~~~~~~~~~~~~~~~~~~~~~~~~~~~~~~~~~~~~~~~~~~~~~}$}
\put(44.5,-1){\Large{$S_1$}}

\put(33,55){$(x,v)$}
\put(33,32){$(x',i)$}
\put(65,45){$(x'',j)$}

\put(29,8){$x^{(i)}$}
\put(61,8){$x^{(j)}$}

\put(29.5,68){$S^{(i)}$}
\put(61.5,68){$S^{(j)}$}

\multiput(30,54)(0.3,0){5}{\line(0,1){1.5}}
\multiput(30,33)(0.3,0){5}{\line(0,1){1.5}}
\multiput(30,54)(0.3,0){5}{\line(0,1){1.5}}
\multiput(30,12)(0.1,0){16}{\line(0,1){1.5}}
\multiput(62,43.5)(0.3,0){5}{\line(0,1){1.5}}
\multiput(30,54)(0.3,0){5}{\line(0,1){1.5}}

\multiput(20.0,12)(10,0){2}{\multiput(0,0)(0.1,0){16}{\line(0,1){1.5}}}
\multiput(52,12)(10,0){2}{\multiput(0,0)(0.1,0){16}{\line(0,1){1.5}}}
\end{picture}
\caption{A schematic representation of the path from $(x,v)\in\bar S$ to $\bar x\in\bar S_1$ in the proof of Corollary \ref{4.3}.} \label{fig:5}
\end{center}
\end{figure}

{\bf Proof. } 
%We have 
%\begin{equation}
%\begin{split}
%x^{(i)}=T_{i_{s(i)}} \dots T_{i_1}(T_i)^{\t-s(i)} x 
%\end{split}
%\end{equation}
%for any $x\in S^{(i)}$ (using of course lazyness of the chain $A$) which means that 
%we can go to a configuration in $\bar S_1$ in a uniform number of time-steps $\t$. 
%This proves the analogue of \eqref{killing},  
%\begin{equation}
%\begin{split}\label{killing1}
%\inf_{(x,v)\in \bar S}M^{\t}((x,v),\bar S_1)\geq \bar \a>0
%\end{split}
%\end{equation}
%with some $\bar \a>0$. 
To prove that, for every joint configuration $\bar x=(x,v)\in \bar S_1$ 
\begin{equation}
\begin{split}\label{communism1}
\inf_{{(x,v)\in \bar S}\atop{\bar x \in {\bar S_1}}}M^{\bar s}((x,v),\bar x )\geq \bar \a'>0
\end{split}
\end{equation} 
we repeat the argument of the independent case with a small modification: 
First we get from $(x,v)$ where $x\in S^{(i)}$ to a point $(x',i)$ 
where $x'\in S^{(i)}$ in $\s(i)$ steps by the second lemma. We warn the reader that the $x'$ might be dependent 
on the particular choice of $x,v$. 
Then we get from $(x',i)$ to $\bar x^{(i)}=(x^{(i)},a^{(i)}_{|V|-1})$ in $2 s(i)$ steps by the first lemma. 
%into $\bar S_1$ with a configuration . 
Then we get from there into the state $j$ of the driving Markov chain by 
means of the connecting string $(c(a^{(i)}_{|V|-1},j),j)$ and adding particles at $j$ 
in an $i$- and $j$-independent number of steps. This is equivalent to saying that 
the driving chain is ergodic. 
Using now lazyness we can go from there into a state
which has a prescribed maximum by adding sufficiently many particles at $j$ (which
will typically be outside of $\bar S_1$) called $(x'',j)$. 
In the third step we go from that 
state into the corresponding state in $\bar S_1$ which has the maximum at $j$ ({see fig. \ref{fig:5} for an illustration of this procedure}). 
Note that these procedures a priori might take a total number of particle droppings 
which could depend on the $(x,v)$. We can produce a number of particle droppings 
$\bar s$ which will do the job for all $(x,v)$ by adding more particles, if necessary,  
at the steps where particles are dropped at the same site. 
This proves the Lemma with 
$\bar s=3 \max_{i\in V}2(s(i) + \s(i))$. 
$\Cox$
\bigskip 

%\begin{figure}[htbp]%put here figure 4 from picture.pdf
   %\centering
  % \includegraphics{} % requires the graphicx package
   %\caption{{\red A schematic representation of the path from $(x,v)\in\bar S$ to $\bar x\in\bar S_1$ in the proof of Corollary \ref{4.3}.}}
   %\label{fig:5}
%\end{figure}

\bigskip

Let us compare to the independent case. Then the $\s(i)$-term is not needed, 
one $s(i)$-term (needed to build up a sufficiently 
high maximum) can be dropped in the independent case, and one $s(i)$-term is just $|V|-1$ (the length 
of a covering string).

The previous considerations given in Section 2 give us now the existence of an invariant 
distribution $\bar \pi$ on $\bar S$, along with the convergence to it, and the bound 
on the coupling matrix 
\begin{equation}
\begin{split} &
\bar D_{t,t'} := \sup_{\bar x,\bar x'\in \bar S} \P_c( \bar X(t') \neq \bar X'(t') | \bar X(t)=\bar x, \bar X'(t)=\bar x')
\end{split}
\end{equation}
where $\P_c$ is the coupling of $\bar X(t)=(X(v),v(t))$ with $\bar X'(t)=(X'(v),v'(t))$. 
We have for times at the distance of the communication time $\bar s$ that 
\begin{equation}
\begin{split} &
\bar D_{t,t+\bar s}\leq 1-(\bar \a')^2 |V|
\end{split}
\end{equation}
and this implies for general times 
\begin{equation}
\begin{split} &
\bar D_{t,t'} \leq (1-(\bar \a')^2 |V|)^{\lfloor\frac{ t'-t }{\bar s}\rfloor}
\end{split}
\end{equation}
The concentration Lemma can be formulated for 
observables $\bar g: \bar S^n \rightarrow \R$ and otherwise stays the same.

All estimates on the maximal height $m_V(t)$ carry over when $\a'$ is replaced by $\bar\a'$ 
and $s$ is replaced by $\bar s$.  This finishes the proof of the Theorem.\bigskip\bigskip
%\newpage 

\section{Extension to layer-dependent particle droppings} 

We will finally give an extension to a model of particle droppings 
which allows also for deposition of particles below the top layer, 
albeit only with a fixed finite depth. This however allows for a large class 
of deposition rules and we will be very general here. 
On the other hand, we want 
to assume a non-nullness condition of particle adsorption at any 
site to the top layer, independently of the configuration and the position 
of the last dropped particle. By the last requirement we exclude part of 
the difficulty dealt with in the case of Markov-chain droppings. 

Take the set $\bar \O$ of finite subsets of $V\times \N$.  
The set $\Phi \in \bar \O$ describes the places where particles are sitting. 
The set $\Phi_v=\{h\in \N: (h,v) \in \Phi\}$ describes the places where particles are sitting 
above the fixed site $v$.  To each $\Phi \in \bar \O$ we associate the height function 
$h(\Phi)= (h_v(\Phi))_{v \in V}$ where $h_v(\Phi)=\max \Phi_v$. 
We introduce the configuration obtained by adding a particle at $i$ applying 
the screening rule by $T_i \Phi =\Phi \cup \{(i, \max_{v: \dist(v, i)\leq 1}h_v(\Phi)+1 ) \}$. 
With this notation we have compatibility with the previously defined action 
on the height profile, i.e. $h(T_i \Phi)= T_i h(\Phi)$. 
Denote the smaller set of configurations obeying nearest neighbor exclusion 
by $\O=\{\Phi \in \O: (v,h)\in \Phi \text{ implies } (w,h)\not\in \Phi \text{ if } w \sim v\}$.
A growth process will be defined on $\O$.  

Let $\Phi(t)$ denote a Markov chain with state space $\O$ and transition 
matrix $M(\Phi,\Phi')$ having the properties \\
1. $\inf_{i\in V, \Phi \in \O}M(\Phi,T_i \Phi)\geq \e>0$ non-null screening rule\\
2. $M(\Phi_1,\Phi')= M(\Phi_2,\Phi')$ if $\Phi_1 \sim_k \Phi_2$ (layer-$k$-depth memory) \\ 
Here we have defined equivalence to the depth $k$, denoted by 
$\Phi_1 \sim_k \Phi_2$ if $h(\Phi_1)=h( \Phi_2)$ (the height profile coincides) and 
$$(\Phi_1)_v \cap [-k + h_v(\Phi_1), h_v(\Phi_1)]= 
(\Phi_2)_v \cap [-k + h_v(\Phi_2), h_v(\Phi_2)]
$$ that is the $k$-depth layer below the height profile coincides. \\
3. $M(\Phi,\Phi')=0$ unless $\Phi'= \Phi\cup \{(v,h)\}$ for a single particle in the $k$-layer below 
the maximum, i.e. $h\in \{-k + h_v(\Phi), (T_v (h(\Phi)))_v\}$.  

To formulate the last condition let us subtract the maximum and define $\Psi_v:=\{x: x+ \max_{w\in V}h_w(\Phi)\in \Phi_v\}$ 
and $\Psi=\cup_{v}(v\times \Psi_v)$
to be the set of occupations shifted by the maximum. 
As a result we have that the height function takes has the maximum zero, \\
i.e. $\max_v h_v(\Psi)=0$. 

Denote by $S$ the set of equivalence classes of images under $\Psi$ w.r.t. looking at the 
$k$-depth layer. So it is the space of possible height-profiles enlarged by the information 
which sites below are occupied, up a depth $k$.

We also want that \\
4. $M(\Phi_1,\Phi'_1)= M(\Phi_2,\Phi'_2)$ if $\Psi(\Phi_1)=\Psi(\Phi_2)$ and 
$\Psi(\Phi'_1)=\Psi(\Phi'_2)$. (height-shift-invariance) 

It is clear that the process has a lift on $S$ as a Markov process.

\begin{theorem} 

\begin{enumerate}
\item The law of Markov process $\Psi(t)$ on the set of $k$-layer depth height-shift equivalence 
classes $S$ converges in total variation to an invariant distribution $\bar \pi$ on $S$.  
\item For each $g:S^n \rightarrow \R$ the random variable 
$g(\Psi(1),\dots,\Psi(n))$ 
obeys the Gaussian concentration bound of Lemma \ref{concentration}, 
with a matrix $D_{t,t'}\leq A e^{-\l (t'-t) }$ for all $t'\geq t$ and zero else.
\item  In particular the function $m_V(t)$ 
obeys the bound \eqref{bound1} and \eqref{bound2} of Theorem 1 for suitable constants $c$ and $C$.  
\end{enumerate}
\end{theorem}

{\bf Outline of Proof. } 
To prove the first assertion of the theorem 
we need to construct a coupling, starting from any two layer configurations $\Psi_1,\Psi_2\in S$. Let us 
do this in several steps. 
Informally speaking one can go first to configurations with the property that the height profile 
takes values in $S_1$ (formulated for the top layer in the same way as we did in the section on 
independent particle droppings) and then
create any desired allowed layer of thickness $k$ by adding only particles which happen to feel the screening which happens
with non-null probability. This can be done for any initial configuration, with the same outcome after sufficiently many steps.
In this way one can produce a coupling between any two initial configurations with a uniform very small probability
$\tilde\alpha>0$ after some very large time $\tilde s$. From that point everything in the proof stays the same.

Now we give some details. Suppose that $h(\Psi)\in S^{i}$ (meaning that the top profile takes the maximum at $i$).   
Applying the sequence of particle additions we look at the resulting configuration 
\begin{equation}
\begin{split}
\tilde \Psi_j:=T_{a^{(j)}_{|V|-1}}\dots T_{a^{(j)}_1}(T_{j})^{|V|-1}T_{a^{(i)}_{|V|-1}}\dots T_{a^{(i)}_1}\Psi 
\end{split}
\end{equation}
By the non-nullness screening condition we know that this has a probability which is bounded uniformly 
below by $\a_1>0$. 
We can be certain that $h(\tilde \Psi)$ is equal to the previously defined $x^{(j)}$ independently 
of the initial condition. However, this might not hold for the $k$-layer below. 
To cure for this we take an arbitrary sequence $a=(v_1,v_2,v_3,\dots ,v_{R})$ in which every vertex appears at least 
$k$ times, and apply the corresponding particle additions using the map $T_{v_j}$. 
This creates a configuration  $\Psi_j= T_{v_{R}}\dots T_{v_1}\tilde\Psi_j$ whose $k$-depth layer is independent 
of the starting configuration $\Psi$. Define now the communication set in layer space by putting 
$\bar S_1= \{  \Psi_j: j\in V\}$. This has the desired properties, and by the previous argument 
proves the first part. 
The second part is a direct application of the concentration statement of Lemma \ref{concentration}. 
To prove the third part we write 
\begin{equation}
\begin{split} &
m_V(t)=t-\sum_{u=1}^{t-1}1_{\max_j h_j(\Phi(u+1))=\max_j h_j(\Phi(u))}=t-\sum_{u=1}^{t-1}1_{A(\Psi(u+1),\Psi(u))}
\end{split}
\end{equation}
where $A(\Psi(u+1),\Psi(u))=\{\#\{j\in V|h(\Psi(u+1))_j - h(\Psi(u))_j \}\leq 1\}$. Note that we have 
written the inequality instead of equality in the last definition in order to account for 
particle depositions below the top layer. 
From here the proof of the concentration of the variable $m_V(t)$ stays the same as in the previous 
two cases. This concludes the proof of the Theorem.  $\Cox$
\bigskip
\bigskip

\subsection*{Acknowledgements}
We would like to thank Stefano Olla and Aernout van Enter for useful discussions.

\end{document}